\documentclass{amsart}
\usepackage{amsmath,amsfonts,amssymb}
\usepackage{graphicx}

\newtheorem{theorem}{Theorem}[section]

\newtheorem{cor}[theorem]{Corollary}

\newtheorem{rem}{Remark}

\theoremstyle{definition}

\theoremstyle{remark}

\numberwithin{equation}{section}

\begin{document}

\title{Weighted Hardy Inequality in $l_2$
}



\author[Ivan Gadjev]{Ivan Gadjev}
 \address{Department of Mathematics and Informatics,
Sofia University,
5 James Bourchier Blvd., 
1164 Sofia, Bulgaria}
\email{gadjev@fmi.uni-sofia.bg}

      

\begin{abstract}
  We study the behaviour of the smallest possible constant $d_n$   in weighted Hardy inequality
$$
\sum_{k=1}^{n}\Big(\frac{1}{k}\sum_{j=1}^{k}a_j\Big)^2 k^\epsilon\le
d(n,\epsilon)\,\sum_{k=1}^{n}{a_k^2}\,k^\epsilon
$$
The exact rate of convergence of $d_n$ is established and 
the ``almost extremal'' sequence is found.
\end{abstract}

\maketitle

\section{Introduction and statement of the results}

In the series of papers  \cite {H1919, H1920, H1925} Hardy proved the next two  inequalities.\\
Suppose that $p>1$, $f(x)\geq0$; that $f$ is integrable over any finite interval $(0,x)$ and $f^p$ is integrable over $(0,\infty)$. Then the next inequality 
\begin{equation}\label{I1}
\int_0^\infty\left(\frac{1}{x}\int_0^xf(t)dt\right)^p\,dx\leq
\left(\frac{p}{p-1}\right)^{p}\,\int_0^\infty f^p(x)dx.
\end{equation}
holds.
This is the original Hardy's integral inequality.

The discrete version of Hardy's inequality reads as
\begin{equation}\label{S1}
\sum_{k=1}^{\infty}\Big(\frac{1}{k}\sum_{j=1}^{k}a_j\Big)^p\leq
\left(\frac{p}{p-1}\right)^{p}\,\sum_{k=1}^{\infty}{a_k^p},\ \ \  \ a_k \geq 0,\quad k\in \mathbb{N}.
\end{equation}

Initially, in his 1920 paper \cite{H1920} Hardy derived  \eqref{S1} with the 
larger constant $p^2/(p-1)$ which appeared on the right-hand side of \eqref{S1}. Later, Landau 
 in the letter \cite{Lan2} dated 1921, published later in \cite{Lan1}, 
was the first to establish  (\eqref{S1}) with the exact constant $(p/(p-1))^{p}$ in the sense that 
there is no smaller one for which \eqref{S1}  holds for every sequence of non-negative numbers $a_k$. 
In the same letter Landau pointed out that equality in \eqref{S1} occurs only for the trivial sequence, that is,
 when $a_k=0$ for every $k \in \mathbb{N}$. 

The lack of nontrivial extremal sequence motivates one to consider the``finite version'' of  \eqref{S1}
\begin{equation}\label{eq16}
\sum_{k=1}^{n}\Big(\frac{1}{k}\sum_{j=1}^{k}a_j\Big)^p\leq
d_n\,\sum_{k=1}^n{a_k^p},\ \ \  \ a_k \geq 0,\quad   k=1,2,...,n.
\end{equation}
 The natural questions are: what are the best constant $d_n$ and the corresponding extremal sequence. 

The behaviour of the constant $d_n$ in \eqref{eq16}  was studied extensively - see, for instance, \cite{HW1}, \cite{HW2},\cite{Wilf}, \cite{Wilf2}, \cite{Pec}, \cite{FS}. 
In \cite{Wilf} Herbert S. Wilf established the exact rate of convergence of
the constant $d_n$ 
\[
d_n=4-\frac{16\pi^2}{\ln^2n}+O\left(\frac{\ln\ln n}{\ln^3n} \right).
\]

In \cite{FS} F. Stampach gave slightly better estimation, i.e.
\[
d_n=4-\frac{16\pi^2}{\ln^2n}+\frac{32\pi^2(\gamma+6\log2)}{\log^3n}+
O\left(\frac{1}{\ln^4n} \right).
\]

In \cite{DIGR} we also studied  the asymptotic behaviour of the constant $d_n$. It was proved there that $d_n$ can be expressed in terms of the smallest zero of a continuous dual Hahn polynomial of degree $n$ (see \cite{Long}), for a specific choice of the parameters, in terms of which these polynomials are defined. Despite that nice interpretation of $d_n$, it was only proved in \cite[Theorem 1.1]{DIGR} that  the next inequalities are true for every natural $n\geq 3$
\begin{equation}\label{eq022}
4\Bigg(1-\frac{4}{\ln n +4}\Bigg)\leq d_n \leq
4\Bigg(1-\frac{8}{(\ln n + 4)^2}\Bigg).
\end{equation}

In all proofs of the above mentioned estimations for the constant $d_n$, the authors substantially used the special properties of the space $l_2$. 
In \cite {DGM} and \cite {GG} we applied a different approach 
which allowed us to give a simpler proof of some of
the mentioned estimations and to find an almost extremal sequence. We proved the next theorem.

\begin{theorem}\label{th2}
Let
\[
a_k=\int_k^{k+1}h(x)dx,
\] 
where 
\begin{equation}\label{eq20}
h(x)=x^{-1/2}\left(2\alpha\cos(\alpha\log x)+\sin(\alpha\log x) \right),\quad 1\le x\le n+1,
\end{equation}
and $\alpha$ is the only solution of the equation
\[
\tan(\alpha\log (n+1))+2\alpha=0\ \ \mathrm{in\ the\ interval}\ \ \left(\frac{\pi}{2\log(n+1)},\frac{\pi}{\log(n+1)} \right).
\]
 Then 
\begin{equation}\label{eq18}
\sum_{k=1}^{n}\Big(\frac{1}{k}\sum_{j=1}^{k}a_j\Big)^2\geq
\frac{4}{1+4\alpha^2}\,\sum_{k=1}^{n}{a_k^2}.
\end{equation}
 \end{theorem}

Since
\[
\frac{4}{1+4\alpha^2}\geq 4-16\alpha^2>4-\frac{16\pi^2}{\log^2(n+1)}
\]
we obtain 
\begin{equation}\label{lb}
d_n\geq 4-\frac{16\pi^2}{\log^2(n+1)}.
\end{equation}
Combining the above mentioned results with the latter we obtain very sharp estimates for $d_n$ for every fixed $n$:
\begin{theorem}\label{Th1.4}
The inequalities 
\begin{equation}\label{estdn}
4-\frac{16\pi^2}{[\log(n+1)]^2} \le d_n \le 
4-\frac{32}{[\log n + 4]^2}
\end{equation}
hold for every natural $n\geq 3$.
\end{theorem}

In \cite {DGM} we also established the exact rate of convergence of $d_n$ and full asymptotic expansion of $d_n$ in terms of the negative powers of $\log n$.
 \begin{theorem}\label{Th1.41}
\begin{equation}\label{asdn}
d_n\sim 
4-\frac{16\, \pi^2}{4\, \pi^2 + (\gamma + 6 \log 2 +\log n)^2}\ \ \mathrm{as}\ \ n\rightarrow \infty,
\end{equation}
where $\gamma$ is the Euler constant, and the following full asymptotic expansion of $d(n,0)$ in terms of the negative powers of $\log n$ holds: For every fixed $m\in \mathbb{N}$, $m \geq 2$,
\begin{equation}\label{asexp}
d_n =  
4- \sum_{k=2}^m c_k\  \frac{1}{[\log n]^k}\, +\, O((\log n)^{-m-1}),\ \ \ n \to \infty, 
\end{equation}
with 
\begin{equation}\label{ck}
c_k = 16 \pi^2  \left(4 \pi^2 + (\gamma+6\log2)^2\right)^{(k-2)/2}\ U_{k-2}\left( -\frac{\gamma+6\log2}{\sqrt{4 \pi^2 + (\gamma+6\log2)^2}} \right),
\end{equation}
where $U_{k-2}$ denotes the Chebyshev polynomial of the second kind of degree $k-2$.
\end{theorem}

The weighted generalised version of  \eqref{eq16} reads as
\begin{equation}\label{S22}
\left(\sum_{k=1}^n \Big|\sum_{j=1}^k a_j \Big|^q u_k\right)^{1/q} \leq
d(n,p,q)\left(\sum_{k=1}^n |a_k|^p\,  v_k \right)^{1/p}
\end{equation}
 where $u_k$, $v_k$ are given weight sequences with positive components.

There are many papers investigating different aspects and applications of Hardy's inequality  \eqref{S22} ;  
see for instance \cite{OK1990}, \cite{KP2003} and the bibliography of \cite{KMP2007}. In most of them
the authors are trying to determine the conditions on the parameters $p$, $q$ and on the weight sequences  $u_k$, $v_k$ under
which the Hardy inequality \eqref{S22} holds for some classes of sequences.

In the present paper we investigate a different aspect of the Hardy's inequality  \eqref{S22}. In the special case when $p=q=2$, $v_k=\epsilon$, $u_k=k^{\epsilon-2}$, $0\le\epsilon<1$ , i.e.  the inequality \eqref{S22} reads as
\begin{equation}\label{S221}
\sum_{k=1}^{n}\Big(\frac{1}{k}\sum_{j=1}^{k}a_j\Big)^2 k^\epsilon\le
d(n,\epsilon)\,\sum_{k=1}^{n}{a_k^2}\,k^\epsilon, \quad 
\end{equation}
The questions we are trying to answer are: what is
the smallest possible constant $d(n,\epsilon)$  for which  \eqref{S221} holds and what is the  extremal sequence for \eqref{S221}.

 In this paper we establish for $0\le\epsilon<1$ the exact rate of convergence of $d(n,\epsilon)$ and find an almost extremal sequence. Our main result is the next theorem.
\begin{theorem}\label{th2}
Let $0\le\epsilon<1$. Then the inequality \eqref{S221} holds, i.e.
\begin{equation*}\label{S21}
\sum_{k=1}^{n}\Big(\frac{1}{k}\sum_{j=1}^{k}a_j\Big)^2 k^\epsilon\le
d(n,\epsilon)\,\sum_{k=1}^{n}{a_k^2}\,k^\epsilon
\end{equation*}
where
\begin{equation}\label{S20}
\frac{4}{(1-\epsilon)^2+4\alpha^2} \le d(n,\epsilon) \le \frac{4}{(1-\epsilon)^2}-
\frac{4\sqrt2}{(1-\epsilon)^2\left(8+(1-\epsilon)^2 \right)}\frac{1}{\ln^2(n+1)} 
\end{equation}
and $\alpha$ is the only solution of the equation
\[
\tan(\alpha\ln (n+1))+\frac{2\alpha}{1-\epsilon}=0\ \ \mathrm{in\ the\ interval}\ \ \left(\frac{\pi}{2\ln(n+1)},\frac{\pi}{\ln(n+1)} \right).
\]
Also, the sequence 
\begin{equation}\label{S21}
a_k=\int_k^{k+1}h_\epsilon(x)dx,
\end{equation}
where 
\begin{equation}\label{S6}
h_\epsilon(x)=x^{-(1+\epsilon)/2}\left(\frac{2\alpha}{1-\epsilon}\cos(\alpha\ln x)+\sin(\alpha\ln x) \right)
\end{equation}
is the "`almost extremal"' sequence, i.e. the next inequality holds
 \begin{equation*}\label{S7}
\sum_{k=1}^{n}\Big(\frac{1}{k}\sum_{j=1}^{k}a_j\Big)^2 k^\epsilon\geq
\frac{4}{(1-\epsilon)^2+4\alpha^2}\,\sum_{k=1}^{n}{a_k^2}k^\epsilon.
\end{equation*}
 \end{theorem}

\begin{rem}
The constant $\frac{4\sqrt2}{(1-\epsilon)^2\left(8+(1-\epsilon)^2 \right)}$ in \eqref{S20}
 is by no means the best one. It could be significantly improved in a lot of ways but that would have made the proof longer and much more complicated. Our goal was to keep the proof  as simple as possible. 
\end{rem}

\begin{cor}\label{cor4}
\begin{equation*}
d(n,\epsilon) \sim 
\frac{4}{(1-\epsilon)^2}-\frac{c}{\ln^2 (n+1)},
\end{equation*}
i.e. there exist constants $c_1=c_1(\epsilon)>0$ and $c_2=c_2(\epsilon)>0$ such that
\begin{equation*}
\frac{4}{(1-\epsilon)^2}-\frac{c_1}{\ln^2 n} \le d(n,\epsilon) \le 
\frac{4}{(1-\epsilon)^2}-\frac{c_2}{\ln^2 n}.
\end{equation*}
\end{cor}
Indeed, it follows from
\begin{align*}
\frac{4}{(1-\epsilon)^2+4\alpha^2}
&=\frac{4}{(1-\epsilon)^2}\frac{1}{1+\frac{4\alpha^2}{(1-\epsilon)^2}}
\geq \frac{4}{(1-\epsilon)^2}\left(1-\frac{4\alpha^2}{(1-\epsilon)^2} \right)\\
&\geq \frac{4}{(1-\epsilon)^2}-\frac{16\pi}{(1-\epsilon)^4\ln^2(n+1)}.
\end{align*}
For the function $h_\epsilon(x)$ defined in \eqref{S6},  we have
\[
0<h_\epsilon(x)<x^{-(1+\epsilon)/2}\alpha\left(\frac{2}{1-\epsilon}+|\ln x| \right)
<\frac{\pi}{\ln(n+1)}x^{-(1+\epsilon)/2}\left(\frac{2}{1-\epsilon}+|\ln x| \right)  
\]
and consequently
\[
a_k\le \frac{\pi}{\ln(n+1)}k^{-(1+\epsilon)/2}\left(\frac{2}{1-\epsilon}+\ln(k+1) \right),\quad 1\le k\le n.
\]
 Then it is obvious, as Landau pointed in his letter to Hardy, that if we let $n\rightarrow\infty$ 
the almost extremal sequence $a_k,\, k=1,2,...$ defined in the Theorem \ref{th2} goes to the zero sequence, i.e.
to the sequence $a_k=0$ for all $k$.

\section{Proof of Theorem \ref{th2}}
Cauchy's inequality yields 
\[
\left(\sum_{i=1}^{k}a_i \right)^2\le\left(\sum_{j=1}^{k}\mu_j^2 \right)
\left(\sum_{i=1}^{k}\frac{a_i^2}{\mu_i^2} \right)
\]
for every pair of  sequences $a_i$ and $\mu_i$, such that $\mu_i\neq 0$ for every $i=1,2,...,n$. 
By multiplying both sides of the latter inequality by $k^{-2+\epsilon}$, summing from 1 to $n$ and
changing the order of summation we get
\begin{align*}
\sum_{k=1}^{n}k^{-2+\epsilon}\left(\sum_{i=1}^{k}a_i \right)^2 \le
\sum_{k=1}^{n}\left[k^{-2+\epsilon}\left(\sum_{j=1}^{k}\mu_j^2 \right) \right]
\left[\sum_{i=1}^{k}\frac{a_i^2}{\mu_i^2} \right]=
\sum_{i=1}^{n}M_ia_i^2i^\epsilon
\end{align*}
where
\[
M_i=i^{-\epsilon}\mu_i^{-2}M_i^*,\quad M_i^*=\sum_{k=i}^{n}k^{-2+\epsilon}\sum_{j=1}^{k}\mu_j^2.
\]
Obviously, 
\[
d(n,\epsilon)\le \max_{1\le i\le n}M_i.
\]
Let us denote for brevity $(1-\epsilon)/2=\beta,\,\,0<\beta\le1/2$. Then
\[
M_i=i^{2\beta-1}\mu_i^{-2}M_i^*,\quad M_i^*=\sum_{k=i}^{n}k^{-(2\beta+1)}\sum_{j=1}^{k}\mu_j^2.
\]
Now we minimize 
\[
\max_{1\le i\le n}M_i
\]
 over the sequences $\mu_i\neq 0,\,\,i=1,...,n$, i.e. find 
\[
\min_{\mu_i\neq 0}\,\max_{1\le i\le n}M_i.
\]
Let
\begin{equation*}\label{S14}
\mu_i=\left(i^{\beta-1}-\frac{1}{A\ln^2(n+1)}\int_i^{i+1}\frac{\ln^2x}{x^{1-\beta}}dx\right)^{1/2}
\end{equation*}
where $A=A(\epsilon)$ is a constant which does not depend on $n$ or $\{a_k\}_1^n$ and will be chosen later.

There exists $\eta_i\in [i,i+1]$ such that
\begin{equation}\label{S15}
i^{1-2\beta}\mu_i^2=i^{-\beta}-\frac{i^{1-2\beta}}{A\ln^2(n+1)}\frac{\ln^2\eta_i}{\eta_i^{1-\beta}}
\geq 
i^{-\beta}-\frac{\ln^2\eta_i}{A\ln^2(n+1)\eta_i^\beta}.
\end{equation}

\begin{align*}
M_i^*=\sum_{k=i}^{n}k^{-(1+2\beta)}
\left(\sum_{j=1}^{k}j^{\beta-1}-\frac{1}{A\ln^2(n+1)}\int_1^{k+1}\frac{\ln^2x}{x^{1-\beta}}dx\right)
=M_1-M_2.
\end{align*}

\begin{align*}
M_1=\sum_{k=i}^{n}k^{-(1+2\beta)}\sum_{j=1}^{k}j^{\beta-1}\le
\sum_{k=i}^{n}k^{-(1+2\beta)}\left(1+\int_1^k x^{\beta-1} \right)
=\frac{1}{\beta}\sum_{k=i}^{n}\frac{k^\beta+\beta-1}{k^{1+2\beta}}
\end{align*}
For $x\geq 1$ the function $g(x)=\frac{x^\beta+\beta-1}{x^{1+2\beta}}$ is decreasing and consequently
\begin{align*}
\beta M_1&=\sum_{k=i}^{n}\frac{k^\beta+\beta-1}{k^{1+2\beta}}\le g(i)+\int_i^ng(x)dx\\
&=\frac{1}{i^{1+\beta}}+\frac{1}{\beta i^\beta}
+(\beta-1)\left(\frac{1}{i^{1+2\beta}}+\frac{1}{2\beta i^{2\beta}} \right)
-\frac{1}{\beta}\left(\frac{1}{n^\beta}-\frac{1-\beta}{2n^{2\beta}} \right).
\end{align*}
But
\begin{align*}
\frac{1}{i^{1+\beta}}
+(\beta-1)\left(\frac{1}{i^{1+2\beta}}+\frac{1}{2\beta i^{2\beta}} \right)
=\frac{1}{i^{1+2\beta}}\left(i^\beta+\beta-1+\frac{\beta-1}{2\beta}i \right)\le 0
\end{align*}
because the function $g_1(x)=x^\beta+\beta-1+\frac{\beta-1}{2\beta}x$ is decreasing for $x\geq 1$.
Then
\begin{equation}\label{S13}
M_1\le \frac{1}{\beta^2 i^{\beta}}-
\frac{1}{\beta^2}\left(\frac{1}{n^\beta}-\frac{1-\beta}{2n^{2\beta}} \right).
\end{equation}

\begin{rem}
It is easy to see that for $\mu_i=i^{(\beta-1)/2}$ we obtain \eqref{S221} with 
$d(n,\epsilon)=\beta^{-2}=\left(\frac{2}{1-\epsilon} \right)^2$.
\end{rem}

\begin{align*}
A\ln^2(n+1)M_2&=\sum_{k=i}^{n}k^{-(1+2\beta)}\int_1^{k+1}\frac{\ln^2x}{x^{1-\beta}}dx\\
&=\beta^{-1}\sum_{k=i}^{n}
\frac{(k+1)^\beta\left(\ln^2(k+1)-2\beta^{-1}\ln(k+1)+2\beta^{-2} \right)-2\beta^{-2}}{k^{1+2\beta}}. 
\end{align*}
The function $x^\beta\left(\ln^2x-2\beta^{-1}\ln x+2\beta^{-2} \right)$ is increasing, so
\begin{align*}
\frac{(k+1)^\beta\left(\ln^2(k+1)-2\beta^{-1}\ln(k+1)+2\beta^{-2} \right)-2\beta^{-2}}{k^{1+2\beta}}\\
\geq \int_k^{k+1}\frac{x^\beta\left(\ln^2x-2\beta^{-1}\ln x+2\beta^{-2} \right)-2\beta^{-2}}{x^{1+2\beta}}
\end{align*}
and
\begin{align}\label{S16}
A\ln^2(n+1)M_2&\geq \beta^{-1}\sum_{k=i}^{n}
 \int_k^{k+1}\frac{x^\beta\left(\ln^2x-2\beta^{-1}\ln x+2\beta^{-2} \right)-2\beta^{-2}}{x^{1+2\beta}}\notag\\
&=\beta^{-1}\int_i^{n+1}\frac{x^\beta\left(\ln^2x-2\beta^{-1}\ln x+2\beta^{-2} \right)-2\beta^{-2}}{x^{1+2\beta}}\notag\\
&\geq \beta^{-1}\int_{\eta_i}^{n+1}\frac{x^\beta\left(\ln^2x-2\beta^{-1}\ln x+2\beta^{-2} \right)-2\beta^{-2}}{x^{1+2\beta}}\notag\\
&=\frac{\ln^2\eta_i+2\beta^{-2}-\beta^{-2}\eta_i^{-\beta}}{\beta^2 \eta_i^\beta}
-\frac{\ln^2(n+1)+2\beta^{-2}-\beta^{-2}(n+1)^{-\beta}}{\beta^2(n+1)^\beta}.
\end{align}
Then from \eqref{S13} and \eqref{S16} we obtain
\begin{align*}
M_i^*\le \frac{1}{\beta^2}\left(\frac{1}{i^\beta}-\frac{\ln^2\eta_i}{A\eta_i^\beta\ln^2(n+1)} \right)
-\frac{1}{\beta^4\eta_i^\beta A\ln^2(n+1)}\left(2-\frac{1}{\eta_i^\beta} \right)-F(n)
\end{align*}
where
\[
F(n)=\frac{1}{\beta^2n^\beta}\left(1-\frac{1-\beta}{2n^\beta} \right)
-\frac{1}{A\beta^2(n+1)^\beta}\left(1+\frac{2}{\beta^2\ln^2(n+1)}-\frac{1}{\beta^2(n+1)^\beta\ln^2(n+1)} \right).
\]
By taking $A=A(\beta)$ big enough we can make $F(n)\geq 0$. 
Indeed, for $n\geq 3$ we have
\begin{align*}
F(n)\geq \frac{1}{\beta^{2}n^{\beta}}
\left[1-\frac{1-\beta}{2n^\beta}-\frac{1}{A}\left(1+\frac{2}{\beta^2\ln^2(n+1)} \right) \right]
\geq \frac{1}{\beta^{2}n^{\beta}}\left(\frac{1}{2}-\frac{1}{A}\left(1+\frac{2}{\beta^2} \right) \right)
\end{align*}
so, it is enough to take $A=2\left(1+\frac{2}{\beta^2} \right)$.

Then
\[
M_i^*\le \frac{1}{\beta^2}\left(\frac{1}{i^\beta}-\frac{\ln^2\eta_i}{A\eta_i^\beta\ln^2(n+1)} \right)
-\frac{1}{\beta^4\eta_i^\beta A\ln^2(n+1)}\left(2-\frac{1}{\eta_i^\beta} \right)
\]
and from \eqref{S15} we obtain
\begin{align*}
M_i\le \frac{1}{\beta^2}-\frac{i^\beta}{\beta^4\eta_i^\beta A\ln^2(n+1)}\left(2-\frac{1}{\eta_i^\beta} \right)\le  \frac{1}{\beta^2}-\frac{1}{2\sqrt2\beta^2(2+\beta^2)\ln^2(n+1)}
\end{align*}
since $\eta_i/i\in[1,2]$. Consequently,
\[
d(n,\epsilon) \le \frac{1}{\beta^2}-\frac{c(\beta)}{\ln^2(n+1)}
=\left(\frac{2}{1-\epsilon} \right)^2-\frac{c(\epsilon)}{\ln^2(n+1)}
\] 
where
\[
c(\epsilon)=\frac{4\sqrt2}{(1-\epsilon)^2\left(8+(1-\epsilon)^2 \right)}.
\]

In order to prove the left inequlity in \eqref{S20} we write the left-hand side of \eqref{S221} in the following way
\[
\sum_{k=1}^{n}\left(\frac{1}{k}\sum_{j=1}^{k}a_j\right)^2 k^\epsilon
=\sum_{i=1}^{n}\left(i^{-\epsilon}a_i^{-1}\sum_{k=i}^{n}\left(k^{\epsilon-2}\sum_{j=1}^{k}a_j \right) \right)
a_i^2i^\epsilon=\sum_{i=1}^{n}M_ia_i^2i^\epsilon
\]
where
\[
M_i=i^{-\epsilon}a_i^{-1}\sum_{k=i}^{n}\left(k^{\epsilon-2}\sum_{j=1}^{k}a_j \right). 
\]
Obviously
\[
\sum_{k=1}^{n}\left(\frac{1}{k}\sum_{j=1}^{k}a_j\right)^2 k^\epsilon
\geq \left(\min_{1\le i\le n}M_i\right)\sum_{k=1}^{n}a_k^2k^\epsilon.
\]
Since the function $h_\epsilon(x)$ is continuous there exists a point $\eta_i\in[i,i+1]$ such that
$a_i=h_\epsilon(\eta_i)$ and consequently
\begin{equation}\label{S23}
i^\epsilon a_i\le \eta_i^\epsilon h_\epsilon(\eta_i)
\end{equation}
For the sequence $a_k$, defined in \eqref{S21} we have
\[
\sum_{j=1}^{k}a_j=\int_1^{k+1}h_\epsilon(x)\,dx=\frac{2}{1-\epsilon}(k+1)^{(1-\epsilon)/2}\sin(\alpha\ln(k+1)).
\]
The function $x^{(1-\epsilon)/2}\sin(\alpha\ln x)$ is increasing so
\begin{align*}
\sum_{k=i}^{n}\left(k^{\epsilon-2}\sum_{j=1}^{k}a_j \right)
&=\frac{2}{1-\epsilon}\sum_{k=i}^{n}k^{\epsilon-2}(k+1)^{(1-\epsilon)/2}\sin(\alpha\ln(k+1))\\
&\geq \frac{2}{1-\epsilon}\int_i^{n+1}x^{(\epsilon-3)/2}\sin(\alpha\ln x)\,dx\\
&\geq \frac{2}{1-\epsilon}\int_{\eta_i}^{n+1}x^{(\epsilon-3)/2}\sin(\alpha\ln x)\,dx
=\frac{4\eta_i^\epsilon h_\epsilon(\eta_i)}{(1-\epsilon)^2+4\alpha^2}.
\end{align*}
From this and \eqref{S23} we obtain
\[
M_i\geq \frac{4}{(1-\epsilon)^2+4\alpha^2},\quad i=1,...,n,\quad\mbox{i.e.}
\quad \min_{1\le i\le n}M_i\geq \frac{4}{(1-\epsilon)^2+4\alpha^2} 
\]
and consequently
\[
d(n,\epsilon)\geq \frac{4}{(1-\epsilon)^2+4\alpha^2}.
\]

The proof of Theorem \ref{th2} is complete.

\section*{Declarations}
The authors have no conflicts of interest to declare that are relevant to the content of this article.

\end{document}